\documentclass[12pt]{article}
\usepackage[english]{babel}
\usepackage{amsmath}
\usepackage{amssymb}
\usepackage{amsfonts}
\usepackage{latexsym,amscd,mathtext}
\title{On an analogue of Schwarz's reflection principle}
\author{V.V. Napalkov (Jr.)}

\begin{document}
\maketitle
\begin{abstract}
We consider the Bergman space on the complex plane.
We prove an analogue of Schwarz's reflection principle for  unbounded quasidisks.
\end{abstract}


\noindent{\it Adress:} Valerii V. Napalkov, Institute of Mathematics with Computer Center
of the Ufa Science Center of the Russian Academy of Sciences
 112, Chernyshevsky str.,Ufa, Russia, 450008.

\noindent{\it e-mail:} vnap@matem.anrb.ru, vnap@mail.ru.

 \newtheorem{lemma}{Lemma}
\newtheorem{theorem}{Theorem}
\newtheorem{definition}{Definition}
\newtheorem{corollary}{Corollary}
\newtheorem{proposition}{Proposition}
\vspace{1cm}

Let $G$   be an unbounded simple connected  Jordan domain in the complex plane and $\partial G$ its  boundary. Assume that $\infty\in\partial G$.  The curve $\partial G$  divides the complex plane into two domains $G$ and $\mathbb C\backslash {\overline G}$. Bergman space $B_2(G)$  consists of all holomorphic functions
$f(z),\, z\in G$ such that
\[
\|f\|^2_{B_2(G)}\stackrel{def}{=}\int_{G}|f(z)|^2\,dv(z)<\infty,
\]
where $dv(z)$ is   the Lebesgue area measure on $G$.
Let us show that the system  $\{\frac{1}{(z-\xi)^2}\}_{\xi\in {\mathbb C\backslash\overline G}}$ belongs to $B_2(G)$  as functions of  the variable $z$.
Suppose that  $\xi_0\in\mathbb C\backslash\overline G$, $d=\frac{\rm dist(\xi_0,\partial G)}{2}$ and
$R_d$ is  a disk with  center at the point $\xi_0$ and radius $d$.
 We have
 \begin{eqnarray}
 \int_{G}\left|\frac1{(z-\xi_0)^2}\right|^2\,dv(z) \le\int_{\mathbb C\backslash {R_d}}\left|\frac1{(z-\xi_0)^2}\right|^2\,dv(z)=\nonumber\\=\int_d^\infty\int_0^{2\pi}\frac1{r^4}\cdot r\, dr d\varphi=
 \frac{4\pi}{d^2}<\infty.
  \end{eqnarray}
   Hence, the function $\frac{1}{(z-\xi_0)^2},\, z\in G$ belongs to $B_2(G)$ and the system  $\{\frac{1}{(z-\xi)^2}\}_{\xi\in {\mathbb C\backslash\overline G}}$ belongs to $B_2(G)$  as functions of  the variable $z$.

  Let us show that the system  $\{\frac{1}{(z-\xi)^2}\}_{\xi\in {\mathbb C\backslash\overline G}}$
 is complete in the space $B_2(G)$.
    By  the Banach Theorem  we must prove that
   the condition
            \[
 \bigl(\tfrac{1}{(z-\xi)^2},g(z)\bigr)_{B_2(G)}=0,\quad \forall \xi\in \mathbb C\backslash\overline G,\,g\in B_2(G),
 \]
implies $g\equiv0$.

Without loss of generality it can be assumed that
  $0\in {\mathbb C}\backslash\overline G$.
The conformal mapping
  $w=\varphi(z)=1/z$ takes domain $G$ ($\mathbb C\backslash{\overline G}$) to bounded domain $G_\varphi$ ($\mathbb C\backslash\overline {G_\varphi}$).
  The mapping $\varphi$ generate an isometry  $T_\varphi$ ( see, e.g., \cite{Gayer})
 \begin{eqnarray}
  f\in B_2(G), \, f(z)\stackrel{T_\varphi}{\longrightarrow} f_\varphi(w)=f(\varphi^{-1}(w))\cdot{\varphi^{-1}}'(w)\in B_2(G_\varphi),\nonumber\\
    (f_\varphi,g_\varphi)_{B_2(G_\varphi)}=(f,g)_{B_2(G)}.\nonumber
 \end{eqnarray}
  The operator $T_\varphi$ takes the system of function
  $\{\frac{1}{(z-\xi)^2}\}_{\xi\in {\mathbb C\backslash\overline G}}\in B_2(G)$ to system
  \[
  \bigl\{\tfrac{1}{w^2}\cdot\tfrac{-1}{(1/w-1/\eta)^2}\bigr\}_{\eta\in
  {\mathbb C }\backslash\overline{G_\varphi}}\in B_2(G_\varphi).
   \]
   We have
   \begin{eqnarray}
 \label{1}
 0=\bigl(\tfrac{1}{(z-\xi)^2},g(z)\bigr)_{B_2(G)}=\bigl(\tfrac{1}{w^2}\cdot\tfrac{1}{(1/w-1/\eta)^2},g_\varphi(w)\bigr)_{B_2(G_\varphi)}=\nonumber\\
  =\eta^2\cdot\bigl(\tfrac{1}{w^2\eta^2}\cdot\tfrac{1}{(1/w-1/\eta)^2},g_\varphi(w)\bigr)_{B_2(G_\varphi)}=\nonumber\\
  =\eta^2\cdot\bigl(\tfrac{1}{(w-\eta)^2},g_\varphi(w)\bigr)_{B_2(G_\varphi)},\, \forall \eta
  \in {{\mathbb C}\backslash\overline G}_\varphi.
   \end{eqnarray}
   It follows from the paper~\cite{NY} that
  $
  \bigl(\tfrac{1}{(w-\eta)^2},g_\varphi(w)\bigr)_{B_2(G_\varphi)},\, \eta\in{\mathbb C\backslash\overline G}_\varphi
   $
is holomorphic function.
Using~(\ref{1}), we get
\begin{equation}
\label{2}
 \bigl(\tfrac{1}{(w-\eta)^2},g_\varphi(w)\bigr)_{B_2(G_\varphi)}=0,\forall \eta\in {{\mathbb C}\backslash\overline G}_\varphi.
\end{equation}
  The system of functions
  $\bigl\{\tfrac{1}{(w-\eta)^2}\bigr\}_{\eta_\in {\mathbb C\backslash\overline G}_\varphi}$ is complete
  in the space $B_2(G)$ (see, e.g.\cite{NY}).
 By ~(\ref{2}), it follows that
  $g_\varphi(w)\equiv0, \,w\in{\mathbb C\backslash\overline G}_\varphi$.
 Hence,
  \[
g(z)\equiv0,\, z\in {\mathbb C\backslash\overline G}.
\] 
 Also, the system  $\{\frac{1}{(z-\xi)^2}\}_{\xi\in {\mathbb C\backslash\overline G}}$
 is complete in the space $B_2(G)$.
 Let  us associate every linear continuous functional $f^*$ on $B_2(G)$ 
 generates by the function $f\in B_2(G)$, to the function
  \[
\widetilde f(\xi)\stackrel{def}{=}\bigl(\tfrac{1}{(z-\xi)^2},f(z)\bigr)_{B_2(G)}=
\int_{\mathbb C\backslash\overline G}\overline{f(z)}\cdot\frac{1}{(z-\xi)^2},\,dv(z) 
\quad \xi\in{\mathbb C\backslash\overline G}.
\]
\begin{definition}
\label{def_HT}
The function $\widetilde f$ is called  Hilbert transform of the functional generated 
by $f\in B_2(G)$.
\end{definition}
  Since the system of functions 
  $\{\tfrac{1}{(z-\xi)^2}\}_{\xi\in {\mathbb C}\backslash\overline G}$  is complete in the space
  $B_2(G)$,  we see that the mapping $f^*\to f$ is injective. The family of functions $\widetilde f$
  forms a space
 \[
 \bigl\{\widetilde f:\, \widetilde f(\xi)=(\tfrac{1}{(z-\xi)^2},f(z))_{B_2(G)}\bigr\}=\widetilde B_2(G),
 \]
 where the induced structure of the Hilbert space is considered, i.e.
 \[
(\widetilde f,\widetilde g)_{\widetilde B_2(G)}
 \stackrel{def}{=}(g,f)_{B_2(G)}, \quad f,g\in B_2(G).
 \]
\begin{definition}
\label{df_quas}
It is said that a bounded simply connected domain $G\subset{\mathbb C}$ is quasidisk if
there exists a constant $C>0$ such that for any $z_1,z_2\in \partial G$
$$
{\rm diam}\,\,l(z_1,z_2)\le C|z_1-z_2|,
$$
where $l(z_1,z_2)$ is the part of $\partial G$ between $z_1$ and $z_2$ which has
the smaller diameter.
\end{definition}

The following theorem is valid (see~\cite{NY})
\begin{theorem}
\label{isomBerg}
\begin{enumerate}
\item{For any $g\in B_2(G)$ the function
$\widetilde g$ lies in $B_2(\mathbb C\backslash\overline G)$ and
 \[
 \|\widetilde g\|_{B_2(\mathbb C\backslash\overline G)}\le \|g\|_{B_2(G)}.
\]}
\item{Suppose the domain $G$ is a quasidisk. Then the Hilbert transform operator acting from $B_2(G)$ to $B_2(\mathbb C\backslash\overline G)$ is surjective operator: for any function
$h\in B_2(\mathbb C\backslash\overline G)$ there exist a unique function
$g\in B_2(G)$ such that $\widetilde g=h$ and
\[
c\|g\|_{B_2(G)}\le\|\widetilde g\|_{B_2(\mathbb C\backslash\overline G)},
\]
were $0<c\le1$ is a some constant.}
\item{ If domain $G$ is not a quasidisk, then the Hilbert transform operator acting from $B_2(G)$ to $B_2(\mathbb C\backslash\overline G)$ is  not surjective operator.
The image of the Hilbert transform operator is dense in 
$B_2(\mathbb C\backslash\overline G)$.}
\end{enumerate}
\end{theorem}

Let the domain $G$ be an unbounded quasidisk and 
$\infty\in\partial G$.

From Ahlfors's theorem (see {\cite {Alf}, p. 48}) it follows that there exist  a quasiconformal reflection $\rho(z)$ such that
\begin{enumerate}
\item{The map $\rho(z)$ is homeomorhphism of the exstended complex plane.}
\item{There exist constants $C_1,C_2>0$ such that
\begin{equation}
\label{s3}
C_1|z_1-z_2|\le|\rho(z_1)-\rho(z_2)|\le C_2|z_1-z_2|,\quad\forall z_1,\,z_2\in \mathbb C.
\end{equation}
}
\item{
\[
\rho(z)=z,\quad z\in \partial G.
\]}
\end{enumerate}
The map $\rho$ is called Ahlfors's quasiconformal reflection.

 Take a point $z_0\in G$.  The functional
$\delta_{z_0}$ is linear continuos functional on $B_2(G)$
(see~\cite{Gayer}), so
$B_2(G)$ is reproducing kernel Hilbert space(\cite{Aron}).
 By $K_{B_2(G)}(z,\xi),\, z,\xi\in G$ we denote  the reproducing kernel of the space $B_2(G)$.(\cite{Aron},\cite{HedKor})

In this paper  we obtain
\begin{theorem}
\label{osn_theor}
Let $G$ be an unbounded quasidisk, $\infty\in\partial G$.
Then there exist  a linear continuous one-to-one operator
$\cal B$ in the space $B_2(G)$ such that
\[
 {\mathcal B}K_{B_2(G)}(z,\rho(\xi))=\tfrac{1}{(z-\xi)^2},\,\forall \xi\in \mathbb C\backslash\overline G.
 \]
\end{theorem}
 \section{Auxiliary information}
\begin{definition}[\cite{Luk}]
Let  $H$ be the Hilbert space over the field ${\mathbb R}$ or ${\mathbb C}$, and $\Omega$
is a space with a countably additive measure $\mu$ (see~\cite{DSH}, p.95--101). The system of elements
$\{e_\omega\}_{\omega\in\Omega}$ is called an orthosimilar system (similar to orthogonal)  with respect to the measure $\mu$ in $H$, if any element $y\in H$ can be represented in the form
\[
y=\int_\Omega (y,e_\omega)_H e
_\omega\,d\mu(\omega).
\]
Here the integral is interpreted as a proper or improper Lebesgue integral of a function with the
values in $H$. In the latter case there is an exhaustion 
$\{\Omega_k\}_{k=1}^\infty$ of the space $\Omega$ possibly
depending on $y$ (it is called suitable for $y$), that the function $(y,e_\omega)_H\cdot e_\omega$ is Lebesgue integrable
on $\Omega_k$ and
\[
y=\int_{\Omega}(y,e_\omega)_H e_\omega\,d\mu(\omega)=\lim_{k\to\infty}(L)\int_{\Omega_k}(y,e_\omega)_H e_\omega\,d\mu(\omega).
\]
Note that all $\Omega_k$ are measurable by $\mu$, 
$\Omega_k\subset\Omega_{k+1}$ for $k\in {\mathbb N}$ and
$\bigcup_{k=1}^\infty\Omega_K=\Omega$.
\end{definition}
{\sc Examples}:
\begin{enumerate}
\item {Any orthonormal basis $\{e_k\}_{k=1}^\infty\subset H$ in an arbitrary Hilbert space H is an orthosimilar
system; any element $y\in H$ can be represented in the form
\[
y=\sum_{k=1}^\infty(y,e_k)e_k.
\]
Here one can take a set $\mathbb N$ as $\Omega$, and as of the measure $\mu$ one can take the
counting measure, i.e. a  measure of the set  from $\mathbb N$ is the amount of different natural numbers
belonging to  the set.}
\item{ Let $H$ be the Hilbert space, $H_1$ is a subspace of $H$, and $P$ is the operator
of orthogonal projection of elements from $H$ onto $H_1$. Let 
$\{e_k\}_{k=1}^\infty\subset H$ be an orthogonal
basis in $H$. Then, the system of elements $\{P(e_k)\}_{k=1}^\infty\subset H_1$
is an orthosimilar system in
$H_1$. (see~\cite{Luk}, Theorem 9). Note, that if $\{e_k\}_{k=1}^\infty$ is an orthogonal basis in $H$, then the
system $\{P(e_k)\}_{k=1}^\infty$, in general, is not an orthogonal basis in $H_1$.}
\item{ Let $H=L_2(\mathbb R)$. The function $\psi\in L_2(\mathbb R)$, $\|\psi\|_{L_2(\mathbb R)}=1$. A system of Morlet wavelets
$\psi_{a,b}(x)=\frac1{\sqrt{ |a|}}\psi\left(\frac{t-b}{a}\right), \quad a\in{\mathbb R}\backslash\{0\},\,b\in {\mathbb R}$ is an orthosimilar system in the space $L_2(\mathbb R)$;
any function $L_2(\mathbb R)$ can be represented in the form
\[
f(x)=\int_{{\mathbb R}\backslash{0}}\int_{\mathbb R}(f(\tau),\psi_{a,b}(\tau))_{L_2(\mathbb R)}\psi_{a,b}(x)\, \frac{d b da}{C_\psi|a|^2},
\]
where $C_\psi>0$ is a constant. The set $\left({\mathbb R}\backslash\{0\}\right)\times{\mathbb R}$ with the measure
$\frac{dbda}{C_\psi|a|^2}$ is taken as the
space $\Omega$ here. (see~\cite{GM},\cite{Luk}).}

\end{enumerate}

\begin{definition}[\cite{Luk}]
 An orthosimilar system is said to be nonnegative if the measure $\mu$ is
nonnegative.
\end{definition}
\begin{theorem}[see~\cite{Nap}]
\label{ThC1}
 Let $H$ be a reproducing kernel  Hilbert space of functions on the domain $G\subset {\mathbb C}$. The norm in the space $H$
has an integral form
\begin{equation}
\label{int}
\|f\|_H=\sqrt{\int_G|f(\xi)|^2\, d\nu(\xi)}
\end{equation}
in the space $H$ if and only if the system of functions $\{K_H(\xi,t)\}_{t\in G}$ is a nonnegative orthosimilar
system with respect to the measure $\nu$ in the space $H$.
\end{theorem}

\begin{lemma}
\label{lem1}
Let a domain $G\subset{\mathbb C}$ be an unbounded quasidisk,
$\infty\in \partial G$.
Then there exist a linear continuous self-adjoint operator 
${\mathcal R}$, which is an automorphism of the
Hilbert space $B_2(G)$, such that
the system 
\[
    \bigl\{{\mathcal R}K_{B_2(G)}(z,\rho(\xi))\bigr\}_{\xi\in {\mathbb C}\backslash\overline G}
    \]
 is an orthosimilar system with  respect to the Lebesgue area measure  
in the space
$B_2(G)$, i.e. any function $f\in B_2(G)$  can be represented in the form
 \[
f(z)= \int_{{\mathbb C}\backslash\overline G} (f(\tau), R_\tau K_{B_2(G)}(\tau,\rho(\xi)))_{B_2(G)}\cdot
    R_zK_{B_2(G)}(z,\rho(\xi))\,dv(\xi),\,z\in G.
      \]

\end{lemma}
{\bf Proof.}
The norm of $B_2(G)$ has the form:
\[
 \|f\|_{B_2(G)}=\sqrt{\int_{G}|f(z)|^2\,dv(z)}.
 \]
After change of variable $z=\rho(\xi)$:
 \[
  \int_{G}|f(z)|^2\,dv(z)=\int_{{\mathbb C}\backslash\overline G}|f(\rho(\xi))|^2\,dv(\rho(\xi)).
  \] 
From relation~\ref{s3}, it follows that
\begin{eqnarray}
  \label{s4}
    C_1^2\int_{{\mathbb C}\backslash\overline G}|f(\rho(\xi))|^2\,dv(z)\le\int_{{\mathbb C}\backslash\overline G}|f(\rho(\xi))|^2\,dv(\rho(z))\le\nonumber\\\le C_2^2\int_{{\mathbb C}\backslash\overline G}|f(\rho(\xi))|^2\,dv(z).
  \end{eqnarray}
Define
 \[
    (f,g)_1\stackrel{def}{=}\int_{{\mathbb C}\backslash\overline G}f(\rho(\xi))\cdot \overline{g(\rho(\xi))}\,dv(\xi); \, \|f\|_1 \stackrel{def}{=}\sqrt{(f,f)_1}.
    \]
Using relation~(\ref{s4}), we get 
 \[
    C_1\|f\|_1\le\|f\|_{B_2(G)}\le C_2\|f\|_{1}, \, \forall f\in {B_2(G)}.
    \]
The norms $\|\cdot\|_{B_2(G)}$,$\|\cdot\|_1$ are equivalent.

From Lemma 1,~\cite{Nap} it follows that there exist
a linear continuous self-adjoint operator 
$T$, 
which is the automorphism of the
Hilbert space $B_2(G)$, such that
\[
  (f,g)_{B_2(G)}=(Tf,g)_{1},\,f,g\in B_2(G).
  \]
We get
  \begin{eqnarray}
  (f,g)_{B_2(G)}=(T f,g)_{1}=\int_{{\mathbb C}\backslash\overline G} Tf(\rho(\xi))\cdot \overline{g(\rho(\xi))}\,dv(\xi)=\nonumber\\=
     \int_{{\mathbb C}\backslash\overline G}(T f (\tau),K_{B_2(G)}(\tau,\rho(\xi)))_{B_2(G)}\cdot
     \overline{(g(\tau),K_{B_2(G)}(\tau,\rho(\xi)))
     _{B_2(G)}}\,dv(\xi).
  \end{eqnarray} 
Consider a  fixed point $z\in G$.
We take the function $g(\tau)=K_{B_2(G)}(\tau,z)$.
Hence, 
\begin{eqnarray}
    f(z)=(f(\tau),K_{B_2(G)}(\tau,z))_{B_2(G)}=(T f(\tau),K_{B_2(G)}(\tau,z))_{1}=\nonumber\\=
     \int_{{\mathbb C}\backslash\overline G}( T f(\tau),K_{B_2(G)}(\tau,\rho(\xi)))_{B_2(G)}\nonumber\times\\
\times\overline{(K_{B_2(G)}(\tau,z),K_{B_2(G)}(\tau,\rho(\xi)))_{B_2(G)}}\,dv(\xi) =\nonumber\\=
    \int_{{\mathbb C}\backslash\overline G} ( T f(\tau),K_{B_2(G)}(\tau,\rho(\xi)))_{B_2(G)}\cdot
     \overline{K_{B_2(G)}(\rho(\xi),z)}\,dv(\xi) =\nonumber\\=
     \int_{{\mathbb C}\backslash\overline G} (T f(\tau),K_{B_2(G)}(\tau,\rho(\xi)))_{B_2(G)}\cdot
     K_{B_2(G)}(z,\rho(\xi))\,dv(\xi).
  \end{eqnarray}
Any function $f\in B_2(G)$ can be represented in the form
\[
    f(z)= \int_{{\mathbb C}\backslash\overline G} 
( T f(\tau),K_{B_2(G)}(\tau,\rho(\xi)))_{B_2(G)}\cdot
     K_{B_2(G)}(z,\rho(\xi))\,dv(\xi).
\]

 Operator $T$ is a self-adjoint operator, so it  has a unique positive square root 
\[
{\mathcal R}:B_2(G)\to B_2(G)
\] 
(see, e.g.,~\cite{Riss},
pp.264, 265) such that $T={\mathcal R}\circ {\mathcal R}$. The operator ${\mathcal R}$ is an automorphism of the space $B_2(G)$ as well.
Therefore,
 \begin{eqnarray}
        \label{s5}
    f(z)= \int_{{\mathbb C}\backslash\overline G} ({\mathcal R}\circ{\mathcal R} f(\tau),K_{B_2(G)}(\tau,\rho(\xi)))_{B_2(G)}\cdot
     K_{B_2(G)}(z,\rho(\xi))\,dv(\xi)=\nonumber\\=
      \int_{{\mathbb C}\backslash\overline G} ({\mathcal R} f(\tau), {\mathcal R}_\tau K_{B_2(G)}(\tau,\rho(\xi)))_{B_2(G)}\cdot
     K_{B_2(G)}(z,\rho(\xi))\,dv(\xi).
   \end{eqnarray}
Since ${\mathcal R}$ is one-to-one operator in the space $B_2(G)$, it follows from~(\ref{s5}) that
 \[
  {\mathcal R}^{-1}f(z)= \int_{{\mathbb C}\backslash\overline G} (f(\tau), R_\tau K_{B_2(G)}(\tau,\rho(\xi)))_{B_2(G)}\cdot
     K_{B_2(G)}(z,\rho(\xi))\,dv(\xi),
   \]
where operator ${\mathcal R}^{-1}$ is the inverse operator to the operator $\mathcal R$.
Using the
theorem from (\cite{DSH}, стр.113), one can demonstrate that
\begin{equation}
 \label{s6}
  f(z)= \int_{{\mathbb C}\backslash\overline G} (f(\tau), R_\tau K_{B_2(G)}(\tau,\rho(\xi)))_{B_2(G)}\cdot
    R_z K_{B_2(G)}(z,\rho(\xi))\,dv(\xi),\,z\in G.
   \end{equation}
Thus, the system 
\[
    \bigl\{{\mathcal R}K_{B_2(G)}(z,\rho(\xi))\bigr\}_{\xi\in {\mathbb C}\backslash\overline G}
    \]
 is an orthosimilar   system with  respect to the Lebesgue area measure  
in the space
$B_2(G)$, i.e. any function $f\in B_2(G)$  can be represented in the form
 \[
f(z)= \int_{{\mathbb C}\backslash\overline G} (f(\tau), R_\tau K_{B_2(G)}(\tau,\rho(\xi)))_{B_2(G)}\cdot
    R_zK_{B_2(G)}(z,\rho(\xi))\,dv(\xi),\,z\in G.
      \]
Lemma~\ref{lem1} is proved.
\begin{lemma}
\label{lem2}
Let a domain $G$ be a quasidisk.
Then there exist a linear continuous self-adjoint operator 
${\mathcal S}$, which is an automorphism of the
Hilbert space $B_2(G)$, such that
the system 
\[
    \bigl\{{\mathcal S}_z\tfrac{1}{(z-\xi)^2})\bigr\}_{\xi\in {\mathbb C}\backslash\overline G}
    \]
 is an orthosimilar   system with respect to the Lebesgue area measure  
in the space
$B_2(G)$, i.e. any function $f\in B_2(G)$  can be represented in the form
 \begin{equation}
 \label{7}
  f(z)= \int_{{\mathbb C}\backslash\overline G} (f(\tau), {\mathcal S}_\tau\tfrac{1}{(\tau-\xi)^2})_{B_2(G)}\cdot
    {\mathcal S}_z\tfrac{1}{(z-\xi)^2} \,dv(\xi),\,z\in G.
   \end{equation}
The reproducing kernel of the space $B_2(G)$ has the form
\begin{equation}
 \label{k7}
  K_{B_2(G)}(z,\eta) =\int_{\mathbb C\backslash\overline G}S_z\tfrac{1}{(z-\xi)^2}\cdot
  \overline{S_\eta\tfrac{1}{(\eta-\xi)^2}} \,dv(\xi), z,\eta\in G.
  \end{equation}
\end{lemma}
{\bf Proof.} 
Since $G$ is a quasisisk, it follows from Theorem~\ref{isomBerg}
that the space $\widetilde B_2(G)$ has  an equivalent norm
 \[
    \|\widetilde g\|_1\stackrel{def}{=}\sqrt{\int_{{\mathbb C}\backslash\overline G}|\widetilde g(\xi)|^2\,dv(\xi)}.
\]
There exist a constant $c>0$ such that
  \[
    c\|\widetilde g\|_{\widetilde B_2({\mathbb C}\backslash\overline G)}\le\|\widetilde g\|_{1}\le\|\widetilde g\|_{\widetilde B_2({\mathbb C}\backslash\overline G)}, \forall \widetilde g\in \widetilde B_2({\mathbb C}\backslash\overline G).
     \]
We use the following theorem from~\cite{Nap}.
\begin{theorem}[\cite{Nap}]
\label{theqv}
In order to introduce into the space $\widetilde B_2(G,\mu)$  a norm
\[
\|\widetilde f\|_\nu=\sqrt{\int_{\mathbb C\backslash \overline G}|\widetilde f(\xi)|^2\,
d\nu(\xi)}
\]
($\nu$ is a nonnegative Borel measure on ${{\mathbb C}\backslash \overline{G}}$),
which is equivalent to the original
one,
 it is necessary and sufficient that there exist a
linear continuous operator $S$, realizing an automorphism of the Banach space $B_2(G,\mu)$ , such that the system 
$\{S\left(\tfrac{1}{(z-\xi)^2}\right)\}_{\xi\in {\mathbb C}\backslash \overline G}$
is an orthosimilar system with  respect to the measure $\nu$ in the space $B_2(G,\mu)$, i.e. any element $f\in B_2(G,\mu)$ can be represented in the form
\[
f(z)=\int_{{\mathbb C}\backslash\overline G}( f(\tau),{ S}_\tau\tfrac{1}{(\tau-\xi)^2})_{B_2(G,\mu)}
{ S}_z \tfrac{1}{(z-\xi)^2}\,d\nu(\xi), \quad z\in G.
\]
\end{theorem}
Let the  measure $\mu$ be the Lebesgue area measure on $G$. 
Let the measure $\nu$ be the Lebesgue area measure on ${\mathbb C}\backslash\overline G$ as measure $\nu$.
By theorem~\ref{theqv}, it follows that the system $\{S_z\tfrac{1}{(z-\xi)^2}\}_{\xi\in {\mathbb C}\backslash \overline G}$
is an orthosimilar system with respect to the Lebesgue area measure on ${\mathbb C}\backslash\overline G$ in the space $B_2(G)$, i.e. any element $f\in B_2(G)$ can be represented in the form
\[
f(z)=\int_{{\mathbb C}\backslash\overline G}( f(\tau),{ S}_\tau\tfrac{1}{(\tau-\xi)^2})_{B_2(G)}
{ S}_z \tfrac{1}{(z-\xi)^2}\,dv(\xi), \quad z\in G.
\]
Fixing a point  $\eta\in G$,
we take the function $f(z)=K_{B_2(G)}(z,\eta)$.
We obtain 
\[
  K_{B_2(G)}(z,\eta) =\int_{\mathbb C\backslash\overline G}S_z\tfrac{1}{(z-\xi)^2}\cdot
  \overline{S_\eta\tfrac{1}{(\eta-\xi)^2}} \,dv(\xi), z,\eta\in G.
  \]
Lemma~\ref{lem2} is proved.

Let us define a linear manifold of functions ${\mathcal L}$ as a set of functions $\varphi\in B_2(G)$ such that there is a
finite set of points $\{\xi_k\}_{k=1}^n\in {\mathbb C}\backslash\overline G$, and a set of complex numbers $\{c_k\}_{k=1}^n\in {\mathbb C}$,
 and the function
$\varphi$ has the form
\[
\varphi(z)=\sum_{k=1}^n c_k R_z K_{B_2(G)}(z,\rho(\xi_k)),\, z\in G.
\]
Thus, ${\mathcal L}$ is a linear span of the system of functions
\[
      \left\{R_z K_{B_2(G)}(z,\rho(\xi))\right\}_{\xi\in {\mathbb C}\backslash \overline G}\subset B_2(G).
\]
\begin{lemma}
\label{lemLn1}
Let $\{\xi_k\}_{k=1}^n$
be a set of n distinct points belonging to the domain 
${\mathbb C}\backslash\overline G$, and $\varphi(z)$ be
a function of the form
\[
\varphi(z)=\sum_{k=1}^n c_k R_z K_{B_2(G)}(z,\rho(\xi_k)),\, z\in G,
\]
where $c_k, k = 1,...,n$ are some constants.
Then, the condition $\varphi(z)\equiv 0$ emplies that 
$c_k=0, k = 1,...,n$.
\end{lemma}
Let us define a linear manifold of functions ${\mathcal M}$ as a set of functions $\psi\in B_2(G)$ such that there is a
finite set of points $\{\xi_k\}_{k=1}^n\in {\mathbb C}\backslash\overline G$, and a set of complex numbers $\{c_k\}_{k=1}^n\in {\mathbb C}$,
 and the function
$\psi$ has the form
\[
\psi(z)=\sum_{k=1}^n c_k S_z\tfrac{1}{(z-\xi_k)^2},
\, z\in G.
\]
Therefore, ${\mathcal M}$ is a linear span of the system of functions
\[
      \left\{S_z \tfrac{1}{(z-\xi)^2}\right\}_{\xi\in {\mathbb C}\backslash \overline G}\subset B_2(G).
\]
\begin{lemma}
\label{lemLn2}
Let $\{\xi_k\}_{k=1}^n$
be a set of n distinct points belonging to the domain 
${\mathbb C}\backslash\overline G$, and $\psi(z)$ be
a function of the form
\[
\psi(z)=\sum_{k=1}^n c_k S_z\tfrac{1}{(z-\xi_k)^2},
\, z\in G,
\]
where $c_k, k = 1,...,n$ are some constants.
Then, the condition $\psi(z)\equiv 0$ emplies that 
$c_k=0, k = 1,...,n$.
\end{lemma}

\section{Proof of Theorem~\ref{osn_theor}.}
Let $G$ be an unbounded quasidisk, $\infty\in \partial G$.
 By Lemma~\ref{lem1},
it follows that the system 
\[
    \bigl\{{\mathcal R}K_{B_2(G)}(z,\rho(\xi))\bigr\}_{\xi\in {\mathbb C}\backslash\overline G}
\]
 is an orthosimilar   system with  respect to the Lebesgue area measure  
in the space
$B_2(G)$, i.e. any function $f\in B_2(G)$  can be represented in the form
 \[
f(z)= \int_{{\mathbb C}\backslash\overline G} (f(\tau), R_\tau K_{B_2(G)}(\tau,\rho(\xi)))_{B_2(G)}\cdot
    R_zK_{B_2(G)}(z,\rho(\xi))\,dv(\xi),\,z\in G.
\]
From Lemma~\ref{lem2} it follows that
the system 
\[
    \bigl\{{\mathcal S}_z\tfrac{1}{(z-\xi)^2})\bigr\}_{\xi\in {\mathbb C}\backslash\overline G}
    \]
 is an orthosimilar   system with respect to the Lebesgue area measure  
in the space
$B_2(G)$, i.e. any function $f\in B_2(G)$  can be represented in the form
\[
  f(z)= \int_{{\mathbb C}\backslash\overline G} (f(\tau), {\mathcal S}_\tau\tfrac{1}{(\tau-\xi)^2})_{B_2(G)}\cdot
    {\mathcal S}_z\tfrac{1}{(z-\xi)^2} \,dv(\xi),\,z\in G.
\]
Define the operator $A:{\mathcal L}\to {\mathcal M}$ by the rule
\[
\varphi(z)=\sum_{k=1}^n c_k R_z K_{B_2(G)}(z,\rho(\xi_k))\to
A\varphi(z)=\sum_{k=1}^n c_k S_z\tfrac{1}{(z-\xi_k)^2},
\, z\in G.
\]
By Lemma~\ref{lem1}, Lemma~\ref{lem2}, it follows that  operator $A$ is well defined.
Thus, operator $A$ is an one-to-one operator  acting from
${\mathcal L}$ to ${\mathcal M}$.
Define the norm on ${\mathcal M}$:
\[
\|A\varphi\|_{\mathcal M}
\stackrel {def}{=}\|\varphi\|_{B_2(G)},\quad \varphi\in {\mathcal L}.
\]
Let $H_{\mathcal M}$ be a complement of ${\mathcal M}$ with respect to the norm $\|\cdot\|_{\mathcal M}$.
$H_{\mathcal M}$ is  a Hilbert space.
Since ${\mathcal L}$ is a dense set in $B_2(G)$(see~\cite{Aron}), we see that the 
operator $A$ has a continuous extension to  Hilbert space $B_2(G)$.
\[
\|Af\|_{H_{\mathcal M}}=\|f\|_{B_2(G)},\quad f\in B_2(G).
\]
Indeed, any function $f\in B_2(G)$ can be  approximated  by the sequence of elements $\{\varphi_k\}_{k\ge 0}$ from ${\mathcal L}$. The image $\{A\varphi_k\}_{k\ge 0}$ is a fundamental sequence in $H_{\mathcal M}$. Then   it follows that 
there exist an element $Af$ in $H_{\mathcal M}$ such that
\[
Af=\lim_{k\to\infty}A\varphi_k.
\]
It is clear that the element $Af$ is well defined.
Therefore, 
the operator
\[
A:B_2(G)\to H_{\mathcal M}
\]
is a linear continuous  bijective operator acting from
$B_2(G)$ to  $H_{\mathcal M}$.
This operator is well defined.
\begin{lemma}
\label{lemrk}
The space $H_{\mathcal M}$ is reproducing kernel Hilbert space.
Reproducing kernel of the space  $H_{\mathcal M}$ has the form:
\[
K_{H_{\mathcal M}}(z,\eta)=\int_{{\mathbb C}\backslash\overline G}S_z\tfrac{1}{(z-\xi)^2}\cdot \overline{S_\eta\tfrac{1}{(\eta-\xi)^2}}\,dv(\xi),\,z,\eta\in G.
 \]
\end{lemma}
{\bf Proof.}
By Lemma~\ref{lem1}, it follows that
any function $f\in B_2(G)$ can be represented in the form:
\[
      f(z)= \int_{{\mathbb C}\backslash\overline G} (f(\tau), R_\tau K_{B_2(G)}(\tau,\rho(\xi)))_{B_2(G)}\cdot
    R_zK_{B_2(G)}(z,\rho(\xi))\,dv(\xi),\,z\in G.
\]
We have
\[
A\circ  R_z K_{B_2(G)}(z,\rho(\xi))= S_z\tfrac{1}{(z-\xi)^2},\quad \xi\in{\mathbb C}\backslash\overline G.
\]
Taking into account theorem (\cite{DSH},p. 113), we obtain
\begin{eqnarray}
\label{10}
 A f(z)=\int_{{\mathbb C}\backslash\overline G} (f(\tau), R_\tau K_{B_2(G)}(\tau,\rho(\xi)))_{B_2(G)}\cdot
    A\circ R_zK_{B_2(G)}(z,\rho(\xi))\,dv(\xi)=\nonumber\\=
 \int_{{\mathbb C}\backslash\overline G} (f(\tau), R_\tau K_{B_2(G)}(\tau,\rho(\xi)))_{B_2(G)}\cdot
 S_z\tfrac{1}{(z-\xi)^2}\,dv(\xi)=\nonumber\\=
\int_{{\mathbb C}\backslash\overline G}(Af(\tau),S_\tau\tfrac{1}{(\tau-\xi)^2})_{H_M}
\cdot S_z\tfrac{1}{(z-\xi)^2},dv(\xi),\quad z\in G.
\end{eqnarray}
Thus, the system
\[
 \bigl\{S_z\tfrac{1}{(z-\xi)^2}\bigr\}_{\xi\in {\mathbb C}\backslash\overline G}
\] 
is an orthosimilar  system with respect to the Lebesgue area measure  
in
$B_2(G)$.
We use the following theorem (see~\cite{Luk},Theorem 1).
\begin{theorem}[\cite{Luk}, An analogue of the Parseval identity] 
Let $\{e_{\omega}\}_{\omega\in \Omega}\subset H$ be a nonnegative
orthosimilar system with respect to the measure $\mu$ in $H$.
Then, for any element $y\in H$ one has
\[
\|y\|_H^2=\int_\Omega|(y,e_\omega)|^2\,d\mu(\omega)
\]
and for any two elements $x,y\in H$ one has
\[
(x,y)_H=\int_\Omega (x,e_\omega)\cdot\overline{(y,e_\omega)}\,d\mu(\omega).
\]
\end{theorem}
This implies that,
\[
\|Af\|_{H_M}^2=\int_{{\mathbb C}\backslash\overline G} |(Af(\tau),S_\tau\tfrac{1}{(\tau-\xi)^2})_{H_M}|^2\,dv(\xi).
\]
Consider the point $z_0\in {\mathbb C}\backslash\overline G$. It follows from Schwarz's inequality that
\begin{eqnarray}
|A_z f(z_0)|\le\sqrt {\int_{{\mathbb C}\backslash\overline G}|(Af(\tau),S_\tau\tfrac{1}{(\tau-\xi)^2})_{H_M}|^2\,dv(\xi)}
\times
\nonumber\\
\times
\sqrt{\int_{{{\mathbb C}\backslash\overline G}}|S_z\tfrac{1}{(z_0-\xi)^2}|^2\,dv(\xi)} =\nonumber\\=
\|Af\|_{H_M}\cdot\sqrt{\int_{{{\mathbb C}\backslash\overline G}}|S_z\tfrac{1}{(z_0-\xi)^2}|^2\,dv(\xi)}.
\end{eqnarray} 
Thus, $H_{\mathcal M}$ is a reproducing kernel Hilbert space.

Using~(\ref{10}), we get
\begin{equation}
\label{rkhm}
K_{H_M}(z,\eta)=\int_{{\mathbb C}\backslash\overline G}S_z\tfrac{1}{(z-\xi)^2}\cdot \overline{S_\eta\tfrac{1}{(\eta-\xi)^2}}\,dv(\xi),\,z,\eta\in G.
 \end{equation}
Lemma~\ref{lemrk} is proved.

By Lemma~\ref{lem2}, 
\[
 K_{B_2(G)}(z,\eta) =\int_{\mathbb C\backslash\overline G}S_z\tfrac{1}{(z-\xi)^2}\cdot
  \overline{S_\eta\tfrac{1}{(\eta-\xi)^2}} \,dv(\xi), z,\eta\in G.
\]
Using~(\ref{rkhm}), we get
\[
K_{H_M}(z,\eta)\equiv K_{B_2(G)}(z,\eta).
\] 
It follows from Moore--Aronszajn's theorem (see~\cite{HedKor}, p.~243) that the space 
$H_{\mathcal M}$ coincides with the space $B_2(G)$.
Thus, the operator $A$ is a linear continuous one-to-one operator
$\cal B$ in the space $B_2(G)$ such that
\[
 A\circ{\mathcal R}K_{B_2(G)}(z,\rho(\xi))={\mathcal S}\tfrac{1}{(z-\xi)^2},\,\forall \xi\in \mathbb C\backslash\overline G.
 \]
Denote
${\mathcal B}\stackrel{def}{=}{\mathcal S}\circ A\circ
{\mathcal R}$.
We have
\[
 {\mathcal B}K_{B_2(G)}(z,\rho(\xi))=\tfrac{1}{(z-\xi)^2},\,\forall \xi\in \mathbb C\backslash\overline G.
 \]
Theorem~\ref{osn_theor} is proved.

We prove the following analogue of Schwarz's reflection principle.
\begin{corollary}
\label{PR}
Let $G$ be an unbounded quasidisk, $\infty\in \partial G$, and  
$\rho$ be Ahlfors's quasiconformal reflection.
If $f\in B_2(G)$ then
$\overline{f(\rho(\xi))}\in B_2({\mathbb C}\backslash\overline G)$.
For any function $g\in B_2({\mathbb C}\backslash\overline G)$ there exists a unique
$f\in B_2(G)$ such that
\[
g(\xi)=\overline{f(\rho(\xi))},\, \xi\in {\mathbb C}\backslash\overline G,
\]
\begin{equation}
\label{nerav1}
C_1\|\widetilde f\|_{B_2({\mathbb C}\backslash\overline G)} \le\|g\|_{B_2({\mathbb C}\backslash\overline G)}\le C_2\|\widetilde f\|_{B_2({\mathbb C}\backslash\overline G)},\, C_1,C_2>0.
 \end{equation}
\end{corollary}
{\bf Proof.}
It follows from theorem~\ref{osn_theor} that there exists  a linear continuous one-to-one operator ${\mathcal B}$ in the space $B_2(D,\mu)$ such that
 \[
{\mathcal B}_z K_{B_2(G)}(z,\rho(\xi))=\tfrac{1}{(z-\xi)^2},\quad \xi\in {\mathbb C}\backslash\overline G
\]
and
 \[
{\mathcal B}^{-1}_z \tfrac{1}{(z-\xi)^2}=K_{B_2(G)}(z,\rho(\xi)),\quad \xi\in {\mathbb C}\backslash\overline G,
\]
where operator ${\mathcal B}^{-1}$ is the inverse operator to the operator $\mathcal B$.
If $f\in B_2(G)$ then
\begin{eqnarray}
 \overline{f(\rho(\xi))} =(K_{B_2(G)}(z,\rho(\xi)),f(z))_{B_2(G)}=({\mathcal B}^{-1}_z\tfrac{1}{(z-\xi)^2}, f(z))_{B_2(G)}=\nonumber\\=(\tfrac{1}{(z-\xi)^2}, {\mathcal B}^{-1}_zf(z))_{B_2(G)}=\widetilde{{\mathcal B}^{-1} f}(\xi), \, \xi\in
{\mathbb C}\backslash\overline G.
\end{eqnarray}
Since $\widetilde{{\mathcal B}^{-1} f}$ belongs to $B_2({\mathbb C}\backslash\overline G)$, we see that
$g(\xi)\stackrel{def}{=}\overline{f(\rho(\xi))}$ belongs 
to $B_2({\mathbb C}\backslash\overline G)$.
By theorem~\ref{isomBerg}, it follows that
\begin{equation}
\label{asy1}
\|g\|_{B_2({\mathbb C}\backslash\overline G)} =\|\widetilde{{\mathcal B}^{-1} f}\|_
{B_2({\mathbb C}\backslash\overline G)}\asymp\|{\mathcal B}^{-1} f\|_{B_2(G)}.
\end{equation}
Symbol $\asymp$ means  that norms are equivalent.
Since 
\begin{equation}
\label{asy2}
 \|{\mathcal B}^{-1} f\|_{B_2(G)}\asymp\|f\|_{B_2(G)},\quad f\in B_2(G)
  \end{equation}
it follows from Theorem~\ref{isomBerg} that
\begin{equation}
\label{asy3}
\|f\|_{B_2(G)}\asymp\|\widetilde f\|_
{B_2({\mathbb C}\backslash\overline G)},\quad f\in B_2(G).
 \end{equation}
The relations ~(\ref{asy1}), (\ref{asy2}), (\ref{asy3}) implies
~(\ref{nerav1}).
Corollary~\ref{PR} is proved.

 \end{document}